\documentclass[10pt]{article}
\usepackage{latexsym}\usepackage{amsbsy}\usepackage{amssymb}

\newtheorem{theorem}{Theorem}[section]

\newtheorem{proposition}[theorem]{Proposition}
\newtheorem{definition}[theorem]{Definition}

\newcommand{\pr}{\noindent{\bf Proof. }}
\newcommand{\ep}{\hspace*{\fill}$\Box$}
\newcommand{\Hom}{\mathrm{Hom}}
\newcommand{\g}{\mathfrak g} 
\newcommand{\TM}{\mathfrak X}
\newcommand{\Real}{\mathbb R}
\newcommand{\R}{\mathbb R}
\newcommand{\N}{\mathbb N} 
\newcommand{\C}{\mathbb C}
\newcommand{\Z}{\mathbb Z} 
\newcommand{\p}{\ensuremath{{\mathcal P}} }
\newcommand{\D}{\ensuremath{{\mathcal D}} } 
\newcommand{\gs}{\ensuremath{{\mathcal G}} }
\newcommand{\ns}{\ensuremath{{\mathcal N}} }
\newcommand{\gR}{{\cal R}}
\newcommand{\G}{{\cal G}}
\newcommand{\CC}{{\cal C}}
\newcommand{\esm}{\ensuremath{{\mathcal E}_M} }
\newcommand{\ad}{\mbox{\rm ad}} 
\newcommand{\id}{\mbox{\rm id}} 
\newcommand{\Ad}{\mbox{\rm ad}}
\newcommand{\cinfty}{{\cal C}^\infty}
\newcommand{\eps}{\varepsilon}
\newcommand{\comp}{\subset\subset} 
\newcommand{\sign}{\mbox{sign}}
\newcommand{\pa}{\partial}
\newcommand{\goesto}{\mapsto}
\newcommand{\Div}{\nabla .}
\renewcommand{\leq}{\leqslant}
\newcommand{\vt}{\left(\frac{\partial}{\partial t}\right)}
\newcommand{\nn}{\nonumber} 
\newcommand{\beq}{ \begin{equation} }
\newcommand{\eeq}{\end{equation} }
\newcommand{\beas}{\begin{eqnarray*}}
\newcommand{\eeas}{\end{eqnarray*}}

\begin{document}

\title{Generalised Connections and Curvature}

\author{Michael Kunzinger, Roland Steinbauer\\
Department of Mathematics, University of Vienna,\\ 
Strudlhofg. \textup{4,  A-1090} Wien, Austria \\ 
e-mail\textup{: \texttt{michael.kunzinger@univie.ac.at}}\\ 
e-mail\textup{: \texttt{roland.steinbauer@univie.ac.at}}
\and\ James A Vickers\\
School of Mathematics, University of Southampton, \\ 
Highfield, Southampton SO\textup{17 1}BJ, United Kingdom \\ 
e-mail\textup{: \texttt{J.A.Vickers@maths.soton.ac.uk}}}
        
\maketitle
                                                                           
\begin{abstract}
The concept of generalised (in the sense of Colombeau) connection on a
principal fibre bundle is introduced. This definition is then used to
extend results concerning the geometry of principal fibre bundles to those
that only have a generalised connection. Some applications to singular
solutions of Yang-Mills theory are given.
\end{abstract}

\section{Introduction}\label{intro}

Recently the theory of Colombeau algebras of generalised functions
\cite{c1,c2} has been applied to a number of areas of geometrical interest
such as the Lie group analysis of partial differential equations (e.g.,
\cite{symm, DKP}) and the study of singular spacetimes in general
relativity (see \cite{vickersESI} for a survey). In order to address these
geometric issues in a satisfactory manner it was necessary to reformulate
the theory of Colombeau algebras to ensure that it was diffeomorphism
invariant.  This was accomplished for the so-called full theory (in which
one has a canonical embedding of Schwartz distributions) in
\cite{found,vim} using calculus in `convenient' \cite{KM}
infinite-dimensional vector spaces.  An alternative approach is to
work instead with the so-called special variant of the theory. Although
this version does not posses a canonical embedding of the space of
distributions it may still be used to model singularities in a non-linear
setting and has the added advantage that it is manifestly diffeomorphism
invariant. As a result there has been considerable use of the special
algebra in addressing geometric problems. See for example chapters 4 and 5
of \cite{book} for applications to symmetries of differential equations and
to general relativity.

In \cite{ndg} the construction of generalised functions in the special
algebra was extended to the theory of generalised sections of vector
bundles. In particular this approach was used to define generalised vector
fields, tensor fields and differential forms thus providing a foundation for
a non-linear theory of distributional geometry.  This work was extended in
\cite{gprg} to give a description of generalised (pseudo-)Riemannian
geometry. In particular the definitions of generalised pseudo-Riemannian
metric, generalised Levi-Civita connection and generalised curvature tensor
in this setting were given.

The aim of the present paper is to give a description of an equally
important area of generalised differential geometry, namely the generalised
theory of connections on principal fibre bundles and on associated vector
bundles. Such a theory turns out to be important mathematically in the
description of characteristic currents for singular connections and
Chern-Weil theory for bundle maps \cite{Harvey} (see also
\cite{Marenich}). In \cite{Harvey} Harvey and Lawson consider connections
on a vector bundle $E$ over a manifold $M$ which are smooth over $M
\setminus \Sigma$ where $\Sigma$ is a closed measure zero subset of $M$
called the singular set. They call such a connection a `singular
connection' and in order to compute its Chern currents the singular
connection $\nabla$ is approximated by a family of smooth connections
$\nabla_\eps$, defined over the whole of $M$, which converge to $\nabla$ on
$M \setminus \Sigma$ as $\eps \to 0$.  The problem with such an approach is
that there are very many different ways of representing $\nabla$ in terms
of a family of smooth functions and it is not always clear how the results
depend on the particular choice. It is precisely this kind of issue that
the Colombeau approach is able to examine using the concepts of equivalence
and association. Singular connections also naturally arise in physics
through singular Yang-Mills fields. The first example of a singular
Yang-Mills field was the `fractionally charged instanton' described by
Forgasc et. al.  \cite{Forgasc1}. This is a (Euclidean) self-dual
Yang-Mills connection on the 4-sphere with a singularity along a
2-sphere. In fact this solution is invariant under a circle action so as
shown by Atiyah \cite{Atiyah} can be viewed as a monopole on hyperbolic
3-space. These singular instanton solutions arise because of the existence
of a locally flat connection with non-trivial holonomy associated with
parallel transport around the singular 2-sphere and are very similar to the
conical singularities studied in \cite{clarke}

The plan of the paper is as follows. In section 2 we recall some basic
facts about the description of both generalised functions and generalised
sections of vector bundles in the special Colombeau algebra.  In section 3
we recall the classical theory of connections, in section 4 we show how
this can be extended to the generalised case and in section 5 we introduce
the concept of the generalised curvature of a connection. Section 6 recalls
the concepts of horizontal lift and holonomy in the classical case and
shows how the theory of generalised functions taking values in a manifold
\cite{gfvm, gfvm2} and of generalised flows \cite{flows} enables these
concepts to be extended to the generalised case. Section 7 introduces the
notion of generalised connection on an associated vector bundle and the
corresponding generalised covariant derivative and demonstrates that in the
case of a generalised linear connection on $TM$ (where $TM$ is
regarded as an associated vector bundle of the frame bundle $LM$) 
this precisely reproduces the definition given
in \cite{gprg}.  Section 8 extends the notion of characteristic class
to the generalised case and in section 9 we end by showing how weakly
singular Yang-Mills connections (used for example to describe the
fractionally charged instanton) may be represented by generalised
connections on the whole space.

\section{Nonlinear distributional geometry}\label{recall}

To begin with, let us fix some notation from differential geometry and briefly
recall the construction of generalised functions in the special
Colombeau algebra. Our principal reference for notation and terminology 
from the theory of algebras of generalized functions is \cite{book}.

In what follows, $M$ will always denote a paracompact,
smooth Hausdorff manifold of dimension $m$ and we denote vector
bundles with base space $M$ by $(E,M,\pi)$ or by $E\to M$, for short. 
The space of smooth sections of a vector bundle $E\to M$ is denoted by
$\Gamma(M,E)$. The
$(r,s)$-tensor bundle over $M$ will be written as $T^r_s(M)$ and we 
denote spaces of tensor fields  by  ${\mathcal
T}^r_s(M):=\Gamma(M,T^r_s(M))$, ${\mathfrak X}(M):=\Gamma(M,TM)$ and ${\mathfrak
X}^*(M):=\Gamma(M,T^*M)$, where $TM$ and $T^*M$ denote the tangent and cotangent
bundle of $M$, respectively, while $\Omega^r(M):=\Gamma(M,\Lambda^r(M))$
is the space of r-forms. For a section $s\in\Gamma(M,E)$ we call
$s^i_\alpha=\Psi^i_\alpha\circ s\circ \psi_\alpha^{-1}$ its $i$-th component 
($1\leq i\leq {n'}$, with $n'$ the dimension of the fibers)
with respect to the vector bundle
chart $(V_\alpha,\Psi_\alpha)$ over the chart $\psi_\alpha$. 

The (special) {\em algebra of Colombeau generalised functions on $M$} 
(\cite{RD, ndg})
is defined as
the quotient algebra $\gs(M):=\esm(M)/\ns(M)$ of the algebra $\esm(M)$ of nets of
smooth functions
$(u_\eps)_{\eps\in(0,1]}\in\CC^\infty(M)^{(0,1]}=:{\mathcal E}(M)$ of
moderate growth modulo the ideal $\ns(M)$ of negligible nets, 
defined, respectively, by the following asymptotic estimates 
\beas
        \esm(M)&:=&\{ (u_\eps)_\eps\in{\mathcal E}(M):\ 
        \forall K\subset\subset M,\ \forall P\in{\mathcal P}(M)\ \exists N\in\N:
  \\&&
        \hphantom{mmmmmmmmmmmmmmimmmmm}\sup_{p\in K}|Pu_\eps(p)|=O(\eps^{-N})\}
   \\[1em]
        \ns(M)&:=& \{ (u_\eps)_\eps\in\esm(M):\, 
        \forall K\subset\subset M,\ \forall m \in\N_0:\
        \sup_{p\in K}|u_\eps(p)|=O(\eps^{m}))\},
\eeas
Here ${\mathcal P}(M)$ is the space of linear differential operators
on $M$. 
Elements of $\gs(M)$ are denoted by $u=[(u_\eps)_\eps]=
(u_\eps)_\eps+\ns(M)$. $\gs(\_)$ is a {\em fine sheaf of differential
algebras} with the operation of Lie derivative (w.r.t.\ smooth vector fields)
defined by $L_\xi u:=[(L_\xi u_\eps)_\eps]$.


$\cinfty(M)$ is embedded into $\gs$ by the ``constant''
embedding $\sigma$, i.e., $\sigma(f):=[(f)_\eps]$, rendering $\CC^\infty(M)$
a subalgebra of $\gs(M)$. There exist injective sheaf
morphisms $\iota: \D'(\_\,)\hookrightarrow \G(\_\,)$ embedding the space
$\D'$ of Schwartz distributions into $\gs$ and which coincide with
$\sigma$ on $\CC^\infty(\_\,)$ (\cite{ndg} Th.\ 2).

The notion of {\em association} is employed to achieve compatibility 
with the respective distributional concepts. A generalised function
$u=[(u_\eps)_\eps]$ is called associated with $0$, $u\approx 0$, if $\int_M u_\eps \mu \to
0$ ($\eps\to 0$) for all compactly supported one-densities $\mu$ and one
(hence every) representative $(u_\eps)_\eps$ of $u$. 
If $\int_M u_\eps \mu \to \langle w,\mu\rangle$ for
some $w\in\D'(M)$ 
then $w$ is called the {\em distributional shadow} (or {\em
macroscopic aspect}) of $u$ and we write $u\approx w$. 

Contrary to the distributional setting, there exists a characterization
of Colombeau generalized functions by their point values on generalized 
points (this property may be viewed as a nonstandard aspect of the
theory). In fact, by inserting 
$p\in M$ into $u\in\gs(M)$ one obtains a well-defined element of ${\mathcal R}$, 
defined as the set of moderate nets of numbers ($(r_\eps)_\eps \in
\R^{(0,1]}$ with $|r_\eps| = O(\eps^{-N})$ for some $N$) modulo negligible
nets ($|r_\eps| = O(\eps^{m})$ for each $m$). More generally, $p\in M$
may be replaced by a generalized point $\tilde p$ and $u, v\in \gs(M)$
coincide if and only if they attain the same point value on each $\tilde p$.
For definitions and details we refer to \cite{point, ndg}.


$\Gamma_\gs(M,E)$, the $\gs(M)$-module of {\em generalised sections}  of a
vector bundle $E\to M$ is defined in a manner similar to that of
generalised functions by  using analogous asymptotic estimates with 
respect to the norm induced by any Riemannian metric on the respective
fibres. Indeed, setting 
$\Gamma_{\mathcal E}(M,E):=(\Gamma(M,E))^{(0,1]}$ and letting 
$\p(M,E)$ denote the space of differential operators $\Gamma(M,E)
\to \Gamma(M,E)$,  we may define $\Gamma_\gs(M,E):=\Gamma_{\esm}(M,E)/\Gamma_\ns(M,E)$,
where
\beas 
\Gamma_{\esm}(M,E)&:=& \{ (s_\eps)_{\eps}\in \Gamma_{\mathcal E}(M,E)
        : \ \forall P\in {\mathcal P}(M,E)\, \forall K\comp M \, \exists
        N\in \N:\\ &&\hspace{4cm}\sup_{p\in K}\|Pu_\eps(p)\| =
        O(\eps^{-N})\}\\ 
\Gamma_\ns(M,E)&:=& \{ (s_\eps)_{\eps}\in
        \Gamma_{\esm}(M,E) : \ \forall K\comp M \, \forall m\in \N:\\
        &&\hspace{3.5cm}\sup_{p\in K}\|u_\eps(p)\| = O(\eps^{m})\}\,
\eeas
For the present purposes the
bundles that we are most interested in are the generalised tensor bundles
$\Gamma_\gs(M,T^r_s(M))$ which we denote $\gs^r_s(M)$ and the generalised
differential r-forms $\Gamma_\gs(M,\Lambda^r(M))$ which we denote
$\Omega^r_\gs(M)$. Generalised sections are denoted by
$s=[(s_\eps)_\eps]=(s_\eps)_\eps+\ns(M,E)$.  

Again, smooth sections of $E\to M$ are embedded as constant nets, i.e.,
via $\Sigma: \Gamma(M,E) \hookrightarrow
\Gamma_\gs(M,E)$ by $\Sigma(s)=[(s)_\eps]$ and we will usually suppress
notationally the embedding $\Sigma$.

$\Gamma_\gs(\_\,,E)$ is a fine sheaf of $\gs(M)$-modules and the
$\gs(M)$-module $\gs(M,E)$ is projective and finitely generated (\cite
{ndg}, Th.\ 5). As $\CC^\infty(M) \subset \gs(M)$,
$\Gamma_\gs(M,E)$ may alternatively be viewed as $\CC^\infty(M)$-module and the two
respective module structures are compatible with respect to the
embeddings. \cite{ndg}, Th.\ 4 provides the following algebraic 
characterisation:
\begin{equation}\label{tensorp}
  \Gamma_\gs(M,E)=\gs(M)\otimes_{\cinfty(M)}\Gamma(M,E)\,,
\end{equation}

Generalised tensor fields may also be viewed as
$\CC^\infty(M)$-multilinear mappings taking smooth vector fields resp.\
one-forms to $\gs(M)$ or as $\gs(M)$-multilinear mappings taking
generalised vector resp.\ covector fields to generalised functions, i.e.,
as $\CC^\infty(M)$- resp.\ $\gs(M)$-modules we have (\cite{ndg}, Th.\ 6)
\beas 
\gs^r_s(M)&\cong&L_{\CC^\infty(M)}({\mathfrak X}^*(M)^r,{\mathfrak
X}(M)^s;\gs(M))\\
\gs^r_s(M)&\cong&L_{\gs(M)}(\gs^0_1(M)^r,\gs^1_0(M)^s;\gs(M)).  
\eeas 
The components of a generalised tensor field $T\in{\gs}^r_s(M)$ 
with respect to the chart $(V_\alpha,\psi_\alpha)$ 
are given as the elements 
\[ 
        T^\alpha\,^{\beta_1\dots \beta_r}_{\gamma_1\dots \gamma_s}
        \,:=\,T|_{V_\alpha}(dx^{\beta_1},\dots
        ,dx^{\beta_r},\pa_{\gamma_1},\dots,\pa_{\gamma_s})
\]
of $\gs(V_\alpha)$. 
As before, compatibility with the purely distributional picture is
accomplished via the notion of association and distributional shadow,
cf.\ \cite{ndg}, Sec.\ 7.

We refer to \cite{ndg,gprg,book} for an introduction to further concepts of 
nonlinear distributional geometry in the Colombeau setting. In particular, 
in section \ref{lifts} we will make use of the space $\gs[M,N]$ of 
Colombeau generalized functions defined on the manifold $M$ and taking 
values in the manifold $N$. For definitions and properties of this space
see \cite{gfvm, gprg, gfvm2}.

\section{The Classical Theory of Connections}\label{classical}

Before developing the concept of a generalised connection we briefly review
the classical theory of connections on a principal fibre bundle following
the approach and (unless stated otherwise explicitly) 
notation of Kobayashi and Nomizu \cite{KN}.

Let $P(M,G)$ be a principal fibre bundle over the manifold $M$ with
structure group $G$ and projection $\pi:P \to M$. Given $x \in M$ then
$\pi^{-1}(x)$ is a closed submanifold of $P$ which is diffeomorphic to $G$
and is called the fibre at $x$. If $p \in P$ and $X \in T_pP$ we say that
$X$ is a {\sl vertical} vector if it is tangent to the fibre through
$p$. The space of vertical vectors at $p$ is denoted $V_p$ and is given by
$V_p=\{X \in T_pP: \pi_*(X)=0 \}$.

We now show how the space of vertical vectors at $p$ is isomorphic to the
Lie algebra $\g$ of $G$. The group $G$ acts freely on $P$ on the right
\begin{equation} \label{action}
{\displaystyle
\begin{array}{rcl}
      R:\,P\times G &\to& P\nn\\ (u,g)&\to& ug=R_gu \,.
\end{array}
}
\end{equation}
Let $A \in \g$ be an element of the Lie algebra of $G$, then $g_t=\exp(tA)$
is a 1-parameter subgroup of $G$ which generates a 1-parameter family of
diffeomorphisms $\phi_t=R_{g_t}$ of $P$. Since the group action moves
points along a fibre, the tangent to the orbit is a vector field $A^*$
which is also tangent to the fibre and is therefore a vertical vector
field. We call $A^*$ the fundamental vector field corresponding to $A$ and
denote the map from $A$ to $A^*$ by $\sigma$,
\begin{equation} \label{fundamental}
{\displaystyle
\begin{array}{rcl}
      \sigma:\,\g&\to& \TM(P)\nn\\ A&\to& A^* \,.
\end{array}
}
\end{equation}
The above map is a Lie algebra homomorphism so that
$\sigma([A,B])=[\sigma(A),\sigma(B)]$ where $[,]$ denotes the bracket in
the Lie algebra $\g$ on the LHS and the Lie bracket between vector fields
on the RHS.

If one restricts the map to a point $p \in P$ then
\begin{equation} \label{iso}
{\displaystyle
\begin{array}{rcl}
      \sigma_p:\,\g&\to& V_p\nn\\ A&\to& A^*(p) \,.
\end{array}
}
\end{equation}
is a linear isomorphism between the Lie algebra and the space of vertical
vectors at $p$.

We will also require the fact that the right action of $G$ on $P$ induces a
corresponding action on the fundamental vector fields which satisfies the
condition \begin{equation} \label{ad}
R_{g*}\left(\sigma(A)\right)=\sigma\left(\ad(g^{-1})A\right) \end{equation} where
$\ad$ denotes the adjoint representation, $\ad: G \to L(\g,\g)$, of $G$ on
$\g$.

We are now in a position to make the following definition.  
\begin{definition}
\label{connection} A connection form $\omega$ is a smooth 1-form on $P$
with values in the Lie algebra $\g$ which satisfies the following
conditions
\begin{itemize}
  \item[(i)] $\forall V \in V_p, \ \omega_p(V)=\sigma^{-1}_p(V)$
  \item[(ii)] $\forall V \in T_{R_gp}P, \
   \omega_{R_gp}(V)=\ad(g^{-1})\omega_p(R^{-1}_{g*}V)$
\end{itemize}
\end{definition}

Given a connection 1-form $\omega$ we may define $H_p$, the horizontal
subspace of $T_pP$, to be the set of vectors annihilated by $\omega$, so
that $H_p=\left\{X \in T_pP: \omega(X)=0\right\}$.  Since
$\omega_p\circ\sigma_p=\id$, the null space of
$\omega_p:T_pP \to \g$ is a $\dim M$ 
dimensional vector space transverse to the fibre. This allows one to
uniquely decompose any vector into a vertical and horizontal part and to
define projections $\pi_v:T_pP \to V_p$ and $\pi_h:T_pP \to H_p$. (Note
that $\pi_v$ depends upon $\omega$ even though the space $V_p$ does not.)
Furthermore (\ref{ad}) and condition (i) above ensure that condition (ii)
holds automatically for vertical vectors. Thus condition (ii) amounts to
the requirement that $R_{g*}$ takes horizontal vectors to horizontal
vectors. This leads to the following alternative definition of a
connection.  
\begin{definition} \label{connection2} A connection on $P$ is an assignment
of subspaces $H_p$ of $T_pP$ such that $p \goesto H_p$ is a smooth 
distribution and which satisfies 
\begin{itemize}
  \item[(i)] $T_pP=V_p\oplus H_p$
  \item[(ii)]$H_{R_gp}=R_{g*}H_p$
\end{itemize}
\end{definition}

\section{The Theory of Generalised Connections}

In order to define a generalised connection we replace the classical
connection 1-form with a generalised 1-form. As we saw in the previous
section a connection form $\omega$ is a 1-form on $P$ with values in the
Lie algebra $\g$. That is 
\begin{equation} 
\omega \in \Omega^1(P,\g)=\Omega^1(P)\otimes_\Real\g \,.  
\end{equation} 
We therefore define a generalised connection 1-form to be an element of 
the space 
\begin{equation}
\Omega^1_\gs(P,\g)=\Omega^1_\gs(P)\otimes_\Real\g \,.  
\end{equation} 
We now turn to the conditions which ensure that the generalised form 
defines a connection.

For a classical connection the first condition is that 
\begin{equation}
\omega_p(V)=\sigma^{-1}_p(V)\in\g, \ \forall V \in V_p \subset T_pP 
\end{equation} 
or equivalently that 
\begin{equation} 
\omega_p(\sigma_p(V^*))=V^*, \ \forall V^* \in \g
\end{equation} 
For the case of a generalised connection we impose this requirement by
demanding that the map $\omega_p\circ\sigma_p$ is the identity in the
appropriate space. More precisely we require 
\begin{equation} 
(p \mapsto \omega_p(\sigma_p(\cdot)))=\id \in \gs(P,L(\g,\g)) \,.  
\end{equation}

We now turn to the second condition which classically takes the form 
\begin{equation}
\omega_{R_gp}(V)=\ad(g^{-1})\omega_p(R^{-1}_{g*}V), \ \forall V \in
T_{R_gp}P 
\end{equation} 
It is less clear how this equation should be interpreted in
the sense of generalised functions. To see this we consider the equivalent
equation at the level of representatives of the generalised connection
1-form 
\begin{equation}
\omega_{R_gp,\eps}(V)=\ad(g^{-1})\omega_{p,\eps}(R^{-1}_{g*}V)+N_\eps 
\end{equation}
where $N_\eps$ is a term which is negligible. The difficulty comes from the
fact that the other terms in the equation depend upon $g$ as well as $p$,
$\eps$ and $V$, so that we cannot simply take $N_\eps$ to be an element of
$\Omega^1_{\mathcal N}(P,\g):=\Gamma_{\mathcal N}(P,T^*P)\otimes_{\Real}\g$
but we must also consider how $N_\eps$ depends upon $g$. The attitude we
will take is to regard $g$ as a parameter and therefore require $N_\eps$ to
satisfy the same estimates as an element of $\Omega^1_{\mathcal N}$ for any
fixed value $g_0$ of $g$, or more generally for $g$ in some compact subset
$L$ of $G$.

More precisely we first define $S$ to be the space of $\g$ valued 1-forms
on $P$ parameterised by elements of $G$ 
\begin{equation} 
S=\{s:G\times P \to T^*P\otimes_{\Real}\g \mid 
p\mapsto s^{g_0}_p \mbox{ is a } \g \mbox{ valued 1-form }  \forall g_0 \in \G \}
\end{equation} 
and then define 
\beas 
S_{\mathcal N}&:=&\{ (N_\eps)_\eps \in S^I :\
\forall K\subset\subset P,\ \forall \D\in\p(P,T^*P)\ \forall L
\subset\subset G\ \forall m\in\N: \\ && \sup_{(p,g)\in K\times L}||\D
N^g_{\eps,p}||=O(\eps^m)\}\,.  
\eeas 
Here $\|\ \|$ denotes the norm induced on the fibres of
$T^*P\otimes_\R \g$ by any Riemannian metric on $P$ and any norm on
$\g$. Clearly $S_{\mathcal N}$ does not depend on these choices.

We are now in a position to give a
definition of a generalised connection.  

\begin{definition} \label{genconnection} A
generalised connection form $\omega$ on a principal fibre bundle $P(M,G)$
is an element of $\Omega^1_\gs(P,\g)$ which satisfies the following
conditions

\begin{itemize}
  \item[(i)] The map $(p \mapsto \omega_p(\sigma_p(\cdot)))=\id \in
\gs(P,L(\g,\g))$.
  \item[(ii)] $\exists (\omega_\eps)_\eps$  with $\omega = [(\omega_\eps)_\eps]$ 
$\exists N \in S_{\mathcal N}$ s.t.\  $\forall g\in G$ $\forall V \!\!\in T_{R_gp}P$:
$$\omega_{R_gp,\eps}(V)=\ad(g^{-1})\omega_{p,\eps}(R^{-1}_{g*}V)+N^g_{p,\eps}(V)\,.$$
\end{itemize}
\end{definition}

Note that $\exists (\omega_\eps)_\eps$ and $\forall (\omega_\eps)_\eps$
is equivalent in (ii) above, i.e., every representative of $\omega$ satisfies
the condition if one does. This is easily verified from the definitions of
$\Omega^1_\gs$ and $S_{\mathcal N}$. 

In formulating this definition we began with the classical definition 
of a connection 1-form on $P$ and adapted the definition in the natural way to
the space of generalised 1-forms on $P$. An alternative but less
satisfactory approach would have been to work with 1-parameter families of
connection 1-forms. We now show that although our definition is more
general it allows one to locally choose representatives
$(\omega_\eps)_\eps$ for $\omega$ such that each $\omega_\eps$ is a
connection 1-form.

\begin{theorem} \label{localrep} 
Let $\omega \in \Omega^1_\gs(P,\g)$ be a generalised
connection 1-form on $P(M,G)$ and let $U$ be a relatively compact open set
$U \subset M$.
Then there exist representatives $(\omega_\eps)_\eps$ of 
$\omega|_{\pi^{-1}(U)}$ 
and an $\eps_0 >0$ such that $\omega_\eps$ is a classical connection 1-form on
$\pi^{-1}(U)$ for all $\eps < \eps_0$.  
\end{theorem}

\pr
Let $f:U \to f(U) \subset P$ be a local section of $P$. Given some
representative $(\omega_\eps)_\eps$ for $\omega \in \Omega_\gs^1(P,\g)$ we
use $\omega_\eps$ to construct a connection 1-form $\tilde\omega_\eps$ on
$\pi^{-1}(U)$. We first define the connection form at points $f(x)$ on the
section and then extend it to other points on the fibre using property (ii)
of a classical connection.

Let $H_{f(x),\eps}=\{X \in T_{f(x)}P: \omega_{f(x),\eps}(X)=0\}$. Then
since $\omega_{f(x),\eps}(\sigma_{f(x)}(A))=A+N_{f(x),\eps}(A)$ for all $A
\in \g$ where $N_{f(x),\eps} \to 0$ as $\eps \to 0$ we see that for
sufficiently small $\eps$ the linear map $\omega_{f(x),\eps}:T_{f(x)}P \to
\g$ is onto. Hence for $\eps$ sufficiently small, the kernel,
$H_{f(x),\eps}$, has dimension $\dim M$ and is transverse to the fibre at
$f(x)$. Since $U$ is relatively compact we can find an $\eps_0 > 0$ so that
this is true for all $x \in U$ and all $\eps < \eps_0$.

We may now use $H_{f(x),\eps}$ to uniquely write any vector $X \in
T_{f(x)}P$ as $X=X_{h \eps}+X_{v \eps}$ where $X_{h \eps}\in H_{f(x),\eps}$ 
and $X_{v \eps} \in V_{f(x)}$. 
This enables us to define the connection form $\tilde
\omega_\eps$ at $f(x)$ (for $\eps < \eps_0$) by 
\begin{equation}
\tilde\omega_{f(x),\eps}(X)=\sigma^{-1}_{f(x)}(X_{v,\eps})\,.  
\end{equation} 

In particular we see that if $X \in H_{f(x),\eps}$ then $\tilde\omega_{f(x),\eps}(X)=0$
and vice-versa, so that the horizontal subspaces of $\omega_\eps$ and
$\tilde\omega_\eps$ agree at $f(x)$,  i.e., 
\begin{equation} 
H_{f(x),\eps}=\tilde H_{f(x),\eps} \,.  
\end{equation} 
Also if $X \in V_{f(x)}$, then
$\tilde\omega_{f(x),\eps}(X)=\sigma^{-1}_{f(x)}(X)$ and hence \begin{equation}
\tilde\omega_{f(x),\eps}(\sigma_{f(x)}(A))=A, \ \forall A \in \g \,.  \end{equation}

We now suppose that $p$ is a general point in $\pi^{-1}(U)$ and let
$x=\pi(p)$. Then we may uniquely define $g(p) \in G$ by the requirement
that $R_{g(p)}f(x)=p$. We now use the transformation property of a
classical connection to define $\tilde\omega_{\eps}$ at $p$ according to
\begin{equation} 
\tilde\omega_{p,\eps}(W):=
\ad(g(p)^{-1})\tilde\omega_{f(x),\eps}(R^{-1}_{g(p)*}W), \ \forall W \in
T_pP 
\end{equation} 
Note in particular that if $W \in V_p$, then $R^{-1}_{g(p)*}(W)
\in V_{f(x)}$ so that 
\beas
\tilde\omega_{p,\eps}(W)&=&\ad(g(p)^{-1})\sigma^{-1}_{f(x)}(R^{-1}_{g(p)*}W)
\\ &=& \ad(g(p)^{-1})\ad(g(p))\sigma^{-1}_p(W) \\ &=& \sigma_p^{-1}(W) \\
\eeas and hence condition (i) for a connection holds for
$\tilde\omega_\eps$ at all points $p \in \pi^{-1}(U)$. 
Since $g(R_hp)=g(p)h$, condition (ii) for a connection follows from
$$
\tilde \omega_{R_hp,\eps}(V) = \ad(h^{-1}g(p)^{-1})\tilde \omega_{f(x),\eps}(R_{g(p)*}^{-1}
R_{h*}^{-1}V) = \ad(h^{-1})\tilde \omega_{p,\eps}(R_{h*}^{-1}V)\,.
$$
Finally, all the maps involved are smooth so
that $\tilde\omega_\eps$ is indeed a connection form on $\pi^{-1}(U)
\subset P$ for each $\eps<\eps_0$.

We next show that $\omega_\eps-\tilde\omega_\eps$ is negligible, so that
$(\tilde\omega_\eps)_\eps$ is also a representative for the generalised
connection $\omega$.

Let $p \in \pi^{-1}(U)$. Since $\tilde\omega_\eps$ is a connection we may
write the tangent space at $P$ as the direct sum 
$T_pP=\tilde H_{p,\eps} \oplus
V_p$ where $\tilde H_{p,\eps}$ is the horizontal subspace of $T_pP$ with 
respect to $\tilde\omega_\eps$.

We first suppose that $X \in V_p$, then 
\begin{equation}
\omega_{p,\eps}(X)=\sigma^{-1}_p(X)+N_{p,\eps}(\sigma^{-1}_p(X)) 
\end{equation} 
and, on the other hand 
\begin{equation} 
\tilde\omega_{p,\eps}(X)=\sigma^{-1}_p(X)\,.
\end{equation} 
Hence if $X \in V_p$ then 
\begin{equation}
\omega_{p,\eps}(X)-\tilde\omega_{p,\eps}(X)=N_{p,\eps}(\sigma^{-1}_p(X)) \,.  
\end{equation} 
Since $N$ is an element of $\ns(P,L(\g,\g))$ and $\sigma$ is
smooth, the difference $\omega_{p,\eps}(X)-\tilde\omega_{p,\eps}(X)$ satisfies
the required estimates for $p$ in the compact subset $K \subset
\pi^{-1}(U)$, when $X \in V_p$, uniformly for $\|X\|\leq 1$ (with
$\|\ \|$ induced on $T_pP$ by any Riemannian metric on $P$).

We now suppose that $X \in \tilde H_{p,\eps}$. We first note that this
implies $R^{-1}_{g(p)*}X \in \tilde H_{f(x),\eps}=H_{f(x),\eps}$ and thus
$\omega_{f(x),\eps}(R^{-1}_{g(p)*}X)=0$. Since $p=R_{g(p)}f(x)$ and
$\omega$ is a generalised connection we may therefore write 
\beas
\omega_{p,\eps}(X) &=& \omega_{R_{g(p)}f(x),\eps}(X) \\ &=&
\ad(g(p)^{-1})\omega_{f(x),\eps}(R^{-1}_{g(p)*}X)+N^{g(p)}_{p,\eps}(X) \\
&=& N^{g(p)}_{p,\eps}(X) \,. \\ 
\eeas 
On the other hand 
\begin{equation}
\tilde\omega_{p,\eps}(X)=0 
\end{equation} 
since $X \in \tilde H_{p,\eps}$. So that \begin{equation}
\omega_{p,\eps}(X)-\tilde\omega_{p,\eps}(X)=N^{g(p)}_{p,\eps}(X) \,.  \end{equation}
Because $p$ ranges over a compact set $K$, $g(p)$ also lies within a
compact set and since $N \in S_{\mathcal N}$, the difference
$\omega_{p,\eps}(X)-\tilde\omega_{p,\eps}(X)$ satisfies the required
estimates when $X \in \tilde H_{p,\eps}$ uniformly for $p$ varying in a
compact set and uniformly for $\|X\|\leq 1$.

Using $\tilde \omega_\eps$, any $X\in T_pP$ with $\|X\|\leq 1$ can be written as $X_v+X_{h\eps}$
with $X_v\in V_p$, $X_{h\eps}\in \tilde H_{p\eps}$ and $\|X_v\|$, $\|X_{h\eps}\|
\leq 1$. Therefore the mapping norm of $\omega - \tilde \omega$ with respect to the norm
induced on the fibres of $TP$ by any Riemannian metric on $P$ satisfies the
$\Omega^1_{\ns}(\pi^{-1}(U),\g)$-estimates, i.e., $\omega = \tilde\omega $ in
$\Omega^1_{\gs}(\pi^{-1}(U),\g)$, as claimed.
\ep

We now show that one may use a generalised connection to define a
generalised projection from $TP$ to the space of vertical vectors.

\begin{definition} 
Let $\omega$ be a generalised connection 1-form on $P$ and let
$(\omega_\eps)_\eps$ be a representation for $\omega$. Then
we define the family of vector bundle homomorphisms 
$\pi_{v,\eps}:TP \to TP$ by 
\beq
\pi_{v,\eps}(V)=\sigma_p(\omega_{p,\eps}(V)), \quad \forall p \in P, \
\forall V \in TpP
\eeq
and set $\pi_{v}:=[(\pi_{v,\eps})_\eps] \in \Hom_\gs[TP,TP]$
\end{definition}
Since $\sigma_p$ is a linear isomorphism which depends smoothly on $p$ and
does not depend on $\eps$, $\pi_v$ is well-defined. Furthermore by
choosing a representation for $\omega$ for which each $\omega_\eps$ is a
connection (which is possible by the previous theorem) we see that
if $V \in T_pP$ then $\pi_{\eps,v}(V) \in V_p$ for all $\eps$ so that $\pi_v$ 
defines a generalised projection onto the space of vertical vectors as claimed.

\begin{definition} 
We define $\pi_h \in \Hom_\gs[TP,TP]$, the projection onto the space of
horizontal vectors, by $\pi_h(V)=\mbox{id}_{TP}-\pi_v$. 
\end{definition}
Since $\pi_v\in \Hom_{\mathop{id}_X}[TP,TP]$, this definition is justified
by the remark following 5.8 in \cite{gfvm2}. Again by taking a
representation for $\omega$ consisting of connections one sees from the
corresponding classical result that $\omega(\pi_h(\cdot))=0$ 
in $\Omega^1(P,\g)$. Moreover, for the distinguished representative 
given by the above theorem each $\pi_{h,\eps}$ projects onto the 
kernel of $\omega_\eps$.

\section{Generalised Curvature} \label{gencurve}

As we have seen at each point $p \in P$ a classical connection $\omega$
defines a projection $\pi_h:T_pP \to H_p$ onto the horizontal subspace at
$p$. This enables one to define the exterior covariant derivative $D$ of a
$\g$ valued  r-form $\phi$ by 
\begin{equation} D\phi(V_1,\dots
,V_{r+1})=d\phi(\pi_hV_1,\dots,\pi_hV_{r+1}) \label{extcov} 
\end{equation} 
Since both $\pi_h$ and $d$ are well defined in the generalised case
the above definition may be immediately extended to generalised forms. 

Classically one applies the exterior
covariant derivative to the connection form to define a 2-form
$\Omega=D\omega$ with values in $\g$, which is called the curvature form of
the connection. An explicit calculation then establishes the structure 
equation 
\begin{equation}
\Omega(U,V)=d\omega(U,V)+[\omega(U),\omega(V)] \label{structure} 
\end{equation}
Note that this differs from the corresponding formula in \cite{KN} due to
our choice of a different convention for exterior product and exterior 
derivative which are those of \cite{CDD}.

We use the above formula to define the curvature of a generalised
connection.

\begin{definition} \label{generalcurve} 
Let $\omega \in \Omega^1_\gs(P,\g)$ be a
generalised curvature 1-form on $P(M,G)$. The generalised curvature two
form $\Omega \in \Omega^2_\gs(P,\g)$ is defined to be 
$[(\Omega_\eps)]$ where 
\begin{equation}
\Omega_{p,\eps}(U,V)=d\omega_{p,\eps}(U,V)+[\omega_{p,\eps}(U),
\omega_{p,\eps}(V)], \ \forall U,V \in T_pP 
\end{equation} 
\end{definition}

An important feature of the classical curvature 2-form is that it
is equivariant under $R_g$
by the adjoint representation. We now show that this remains
true for a generalised connection.

\begin{theorem} \label{adinv} 
Let $\Omega$ be the generalised curvature of a
generalised connection $\omega$ then \begin{equation}
(R^*_g\Omega)(U,V)=\ad(g^{-1})\Omega(U,V) \end{equation} where $R^*_g$ is the
pull-back of the generalised 2-form $\Omega$ by the map $R_g:P \to P$ and
is defined by $R^*_g\Omega:=[(R^*_g\Omega_\eps)]$.  
\end{theorem} 
Note that this
equation is to be interpreted in the same generalised sense as the
corresponding equation for the connection. Namely, given a representative
$(\Omega_\eps)_\eps$ for $\Omega$ and a fixed value of $g \in G$ then the
representatives for the left and right hand sides differ by an element of
$\Omega^2_{\mathcal N}(P,\g)$

\pr

Let $U,V \in T_pP$ then 
\beas 
(R^*_g\Omega)_{p,\eps}(U,V) &=&\Omega_{R_gp,\eps}(R_{g*}U,R_{g*}V) \\ 
&=&(d\omega)_{R_gp,\eps}(R_{g*}U,R_{g*}V)+[\omega_{R_gp,\eps}(R_{g*}U),
\omega_{R_gp,\eps}(R_{g*}V)] \\ 
&=&(R^*_g(d\omega))_{p,\eps}(U,V)+[\omega_{R_gp,\eps}(R_{g*}U),
\omega_{R_gp,\eps}(R_{g*}V)] \\ 
&=&(d(R^*_g\omega))_{p,\eps}(U,V)+[\omega_{R_gp,\eps}(R_{g*}U),
\omega_{R_gp,\eps}(R_{g*}V)] \\ 
&=&\ad(g^{-1})(d\omega)_{p,\eps}(U,V)+dN^g_{p,\eps}(R_{g*}U,R_{g*}V) \\
&&+[\ad(g^{-1})\omega_{p,\eps}(U)+N^g_{p,\eps}(R_{g*}U),\ad(g^{-1})
\omega_{p,\eps}(V)+N^g_{p,\eps}(R_{g*}V)] \\ 
&=&\ad(g^{-1})\{d\omega_{p,\eps}(U,V)
+[\omega_{p,\eps}(U),\omega_{p,\eps}(V)]\} \\
&&+dN^g_{p,\eps}(R_{g*}U,R_{g*}V) +
[N^g_{p,\eps}(R_{g*}U),\ad(g^{-1})\omega_{p,\eps}(V)]\\
&&+[\ad(g^{-1})\omega_{p,\eps}(U),N^g_{p,\eps}(R_{g*}V)] \\
&=&\ad(g^{-1})\Omega_{p,\eps}(U,V)+M^g_{p,\eps}(U,V) \,, 
\eeas 
where for fixed $g \in G$, $M^g(U,V) \in \Omega^2_{\mathcal N}(P,\g)$ since
for fixed $g$, $N^g \in \Omega^1_{\mathcal N}(P,\g)$ which is a
differential ideal.  
\ep

The second key property of the curvature $\Omega$ is that it is a {\sl
horizontal form}. That is $U \in V_p \Rightarrow \Omega_p(U,V)=0, \ \forall
V \in T_pP$. This follows immediately from the definition using the
exterior covariant derivative but is less obvious from the definition using
the structure equation. We now show that this result remains true in the
generalised case. Note that even in the generalised case the definition of
the vertical subspace $V_p$ does not depend upon $\eps$ since it is defined
using $\pi:P \to M$ and does not depend upon the generalised connection,
(of course this is not true for the horizontal subspace).

\begin{theorem} \label{horizontal} 
Let $\Omega \in \Omega^2_\gs(P,\g)$ be a generalised
2-form, then $\Omega$ is a horizontal form.  
\end{theorem}

\pr 
By theorem \ref{localrep} we may take local representatives
$\omega_\eps$ of the generalised connection $\omega$ which are themselves
connection forms. The result then follows from the classical result.
\ep

The significance of the above two results is that the curvature is a {\sl
tensorial} 2-form of type $\ad$, (unlike the connection which is only
pseudo-tensorial). This means that the pull-back of $\Omega$ onto $M$ by a
local section transforms under the adjoint action on a change of section
and does not have any inhomogeneous terms (unlike the connection
$\omega$). This will be important in section \ref{characteristic} where we
define characteristic classes for generalised connections. Another
important result that we will use later is Bianchi's identity for the
generalised curvature $\Omega$.

\begin{theorem}[Bianchi's identity]  \label{Bianchi}
$D\Omega=0$ 
\end{theorem}

\pr
By theorem \ref{localrep} we may take local representatives
$\omega_\eps$ of the generalised connection $\omega$ which are themselves
connection forms. The result then follows from taking the exterior
derivative of the structure equation (\ref{structure}) and using the fact
that by definition the connection vanishes on horizontal vectors.
\ep

\section{Horizontal Lifts and Holonomy} \label{lifts}

Given a connection form $\omega$ on $P(M,G)$ we now show how this allows
one to lift a curve $\gamma$ on the base manifold $M$ to a horizontal curve
$\bar\gamma$ on the bundle $P$.  
\begin{definition} \label{lift} 
Let $\gamma:[a,b] \to M$,
$t \mapsto \gamma(t)$ be a smooth curve on $M$.  A horizontal lift of $\gamma$
is a curve $\bar\gamma:[a,b] \to P$ such that
\begin{itemize}
\item[(i)] $\pi\circ\bar\gamma=\gamma$,
\item[(ii)]${\displaystyle \omega\left(\bar\gamma_*\vt\right)=0}$.
\end{itemize}
\end{definition}

Given some point $p_0 \in \pi^{-1}(\gamma(a))$ there exists a unique
lift $\bar\gamma$ such that $\bar\gamma(a)=p_0$ (Proposition 3.1 of
Chapter II in \cite{KN}). We first choose a reference curve $f:[a,b] \to M$
satisfying $f(a)=p_0$ and $\pi\circ f=\gamma$, and look for a curve $g(t)$ 
in $G$ such that $\bar\gamma(t)=R_{g(t)}f(t)$ is horizontal. This is the
case if the curve $g(t)$ in $G$ satisfies 
\begin{equation} \label{de} 
\dot g(t)g(t)^{-1}=-\omega\left(f_*\vt\right) \,.  
\end{equation} 
Thus finding a horizontal lift for $\gamma$ which starts at $p_0$ is 
equivalent to obtaining a solution to (\ref{de}) on $[a,b]$ which
satisfies the initial condition $g(a)=e$. Since the right hand side is
smooth and $[a,b]$ is compact there exists a unique solution to this
equation, and hence there exists a unique lift $\bar\gamma$ starting from
$p_0$. Note that if one chooses a different initial point $p_1=R_{g_0}p_0$
then the corresponding lift is $R_{g_0}\bar\gamma$.

The concept of horizontal lift is used to define the holonomy associated
with a closed curve on $M$. Let $\gamma: [a,b] \to M$ be a closed curve
with $\gamma(a)=\gamma(b)=m$. If $\bar\gamma: [a,b] \to P$ is some
horizontal lift of $\gamma$ then $\bar\gamma(a)$ and $\bar\gamma(b)$ will
be points in the same fibre, so that they define an element $g$ of $G$ such
that $\bar\gamma(b)=R_g\bar\gamma(a)$. We call $g$ the element of holonomy
generated by $\bar\gamma$.  Changing the starting point $\bar\gamma(a)$ by
an element of $G$ simply translates the entire lift, so we see that the
holonomy does not depend upon which horizontal lift is taken but is
determined by $\gamma$. Furthermore the holonomy does not depend upon the
particular parameterisation of $\gamma$, so the holonomy only depend upon
the loop and the connection $\omega$. If one traverses the loop in the
opposite direction the holonomy generated is the inverse of $g$. Similarly
if loops $\gamma_1$ and $\gamma_2$ are two loops based at $m$ which
generate elements of holonomy $g_1$ and $g_2$ respectively then the loop
consisting of $\gamma_1$ followed by $\gamma_2$ generates $g_2g_1$. This
leads to the important notion of holonomy group which encodes considerable
information about the curvature of $P$.

The above concepts can be extended to the generalised setting. We start
with the concept of horizontal lift 
\begin{definition} \label{genlift} 
Let $\gamma \in \gs[[a,b],M]$ 
be a generalised curve in $M$. A generalised
horizontal lift of $\gamma$ is a curve $\bar\gamma \in \gs[[a,b],P]$ such
that

\begin{itemize}
\item[(i)] $\pi\circ\bar\gamma =\gamma$ in $\gs[[a,b],M]$.
\item[(ii)]${\displaystyle
\omega_{\bar\gamma(.)}\left(\bar\gamma_{*}\vt\right)=0}$  in $\gs[[a,b],\g]$.
\end{itemize}
\end{definition}
We note that $\omega_{\bar\gamma}$ is an element of the hybrid space
$\gs^h([a,b],T^*P\otimes \g)$, cf.\ \cite{gprg}, Def.\ 4.4 and 
Th.\ 4.5 (i). On the level of representatives, the above conditions
translate into
\begin{itemize}
\item[(i)] $\pi\circ\bar\gamma_\eps=\gamma_\eps+r_\eps$, where $r_\eps \in
  {\mathcal N}([a,b], M)$
\item[(ii)]${\displaystyle
\omega_{\bar \gamma_\eps(t)}\left(\bar\gamma_{*,\eps}\vt\right)=N_{\eps}(t)}$, where
$(N_\eps)_\eps \in {\mathcal N}([a,b], \g)$.
\end{itemize}

In order to show that generalised lifts exist we may adapt the classical
proof and show that there exist generalised flows $g_\eps(t)$ such that
\begin{equation}
\dot g_\eps(t)g_\eps(t)^{-1}=-\omega_\eps\left(f_*\vt\right) \quad \hbox{in
  $\G$} \,.  
\end{equation}
Again this equation must be interpreted in the generalised sense described
in \cite{flows}.

These ideas will be developed further in a future paper. For the moment we
simply remark that in order to describe the generalised geometry of fibre
bundles one requires a theory which not only permits the multiplication of
generalised functions (in order to define the curvature), but can also deal
with generalised functions (such as $\bar\gamma$) valued in a
manifold, the composition of generalised functions (such as
$\omega\circ\bar\gamma_{*}$) and generalised flows on manifolds
(such as the equation for $g(t)$) all of which are available in the
Colombeau theory as described in \cite{ndg}, \cite{gprg}, \cite{gfvm},
\cite{gfvm2} and \cite{flows}.
 
\section{Generalised Connections in Vector Bundles and Covariant Derivatives} \label{vector}


In many applications one is interested in a connection or covariant
derivative defined on an associated vector bundle rather than a connection
on a principal bundle. For example in gauge theory the matter fields or
Higgs fields are defined on associated bundles and in general relativity
the connection used is defined on the tangent bundle $TM$, which may be
regarded as an associated bundle of the frame bundle $LM$. We start by
defining a covariant derivative (or a connection) on a vector bundle $E$.

A covariant derivative on a vector bundle $E$ is a map \beas \TM(M) \times
\Gamma(M,E) &\to& \Gamma(M,E) \\ (X,V) &\mapsto& \nabla_XV \\ \eeas which
satisfies the conditions
\begin{itemize}
\item[(1)] $\nabla_XV$ is $\Real$-linear in $V$,

\item[(2)] $\nabla_XV$ is $C^\infty(M)$-linear in $X$,

\item[(3)] $\nabla_X(\lambda V)=\lambda\nabla_XV+X(\lambda)\cdot V, \quad
\forall \lambda \in C^\infty(M)$.
\end{itemize}
We now briefly consider the way in which a connection on a principal bundle
may be used to define a covariant derivative on an associated vector bundle
and show how this construction may be extended to the case of a generalised
connection.

We begin by reviewing the construction of an associated vector bundle.  Let
$P(M,G)$ be a principal fibre bundle and let $\rho: G \to GL(n, \Real)$, $g
\mapsto \rho_g$ be a representation of $G$ on $\Real^n$.  Then on the
product manifold $P \times \Real^n$ we may define a right action of $G$
according to \beas G \times (P \times \Real^n) &\to & P \times \Real^n \\
(g,(u,\xi)) &\mapsto& (R_gu, \rho^{-1}_g\xi) \\ \eeas Now let $E=P \times_G
\Real^n$ be the quotient space under this action and let $u\xi \in E$
represent the equivalence class of $(u,\xi) \in P \times \Real^n$. The
projection which takes $(u,\xi) \in P \times G \mapsto \pi(u) \in M$
induces a projection $\pi_E:E \to M$. Every point $x \in M$ has a
neighbourhood $U$ such that $\pi^{-1}(U)$ is isomorphic to $U \times
G$. Indeed let $\Phi:\pi^{-1}(U) \to U \times G$ be such a local
trivialisation then we can use this to define a local trivialisation of $E$
through a map $\Psi:\pi^{-1}(U) \to U \times \Real^n$ which takes $u\xi$ to
$\rho(\Phi_2(u))\xi$. Note that this does not depend upon the choice of
representation of the element of $E$ since an equivalent representation
$(R_hu)(\rho^{-1}_h\xi)$ gets mapped to the same point. We now introduce a
differentiable structure on $E$ by requiring that $\pi^{-1}_E(U)$ is
actually an open submanifold and $\Psi$ is a diffeomorphism. We then say
that $E(M,\Real^n,G,P,\rho)$ is a vector bundle associated to the principal
bundle $P(M,G)$.

We next briefly recall how a (classical) connection on $P$ may be used to
define a covariant derivative on $E$. Let $\Phi: \pi^{-1}(U)\to U \times G$
and $\Psi:\pi_E^{-1}(U) \to U \times \Real^n$ be the local trivialisations
described above. Let $\gamma$ be a curve connecting $\gamma(0)=x_0$ and
$\gamma(1)=x_1$ in $U$ and $\bar \gamma$ be a horizontal lift of the curve
to $P$. Then we may use $\gamma$ to define an isomorphism $\tau ^0_1$
between $\pi_E^{-1}(x_0)$ and $\pi_E^{-1}(x_1)$. Let $V_0 \in
\pi_E^{-1}(x_0)$ then we define $V_1=\tau^0_1(V_0)$ by requiring that \begin{equation}
\xi_1=\rho(g_1g^{-1}_0)\xi_0 \end{equation} where $\Psi(V_i)=(x_i, \xi_i)$ and
$\Phi(\bar \gamma(i))=(x_i, g_i)$ for $i=1,2$.

It is a classical result that $V_1$ does not depend upon either the choice
of trivialisation or upon the particular lift $\bar \gamma$ chosen but only
on the connection on $P$ and the curve $\gamma$. More generally given a
curve $\gamma$ defined on some interval $J$ we may define $\tau^s_t(V)$ to
be the result of parallely propagating the vector $V$ from $\gamma(s)$ to
$\gamma(t)$. We use this to define the covariant derivative of a field
$V(x) \in \Gamma(M,E)$ in the direction of the tangent to the curve at the
point $\gamma(0)$ as \begin{equation} \nabla_{\dot \gamma(0)}V(\gamma(0))=\lim_{h \to
0}{1 \over h}\left[\tau^{h}_0(V(\gamma(h)))-V(\gamma(0))\right] \end{equation} It is
not hard to see that this only depends upon the direction of the tangent
$X=\dot \gamma(0)$ and not the curve $\gamma$ so we may use the above to
define $\nabla_XV$ at $\gamma(0)$. By using this formula at every point $x$
we may define the covariant derivative of the field $V \in \Gamma(M,E)$
with respect to the vector field $X \in \TM(M)$. It is readily verified
that the covariant derivative defined in this way satisfies properties
(1)--(3).

Now let $(V_\alpha, \Psi_\alpha)$ be a vector bundle chart for $E$ with
coordinates $x^i$, $i=1,\dots, m$ on $M$ and let  $e_A$, $A=1,\dots n$ be
the fields on $E$ induced by the canonical basis on $\Real^n$. Then we may
define the $mn^2$ functions  $\Gamma^B_{iA}$ 
(which are the coefficients of the connection on
$E$) by \begin{equation} \nabla_{\partial_i}e_A=\Gamma^B_{iA}e_B \end{equation} An explicit
formula for the connection coefficients in terms of the connection 1-form
on $P$ may also be given. Let $x_0 \in V_\alpha$ and let $\Phi:\pi^{-1}(U)
\to U \times G$ be some local trivialisation of $P$. Choose $p_0 \in
\pi^{-1}(x_0)$ so that $\Phi(p_0)=(x_0,e)$ and define the section $s:
V_\alpha \to \pi^{-1}(V_\alpha)$ by $s(x)=\Phi^{-1}(x,e)$;
then at $x=x_0$
the connection coefficients $\Gamma^B_{iA}$ are given by \begin{equation}
\Gamma^B_{iA}e_B=((D_e\rho)(s^*\omega(\partial_i)))(e_A) \label{coeffs}
\end{equation} The covariant derivative of a field $V$ is then given by
\begin{equation}
\nabla_XV=\left(X^iV^B_{,i}+\Gamma^B_{iA}V^AX^i\right)e_B \label{deriv}
\end{equation}
It is not hard to show that a connection defined by (\ref{coeffs}) and
(\ref{deriv}) does not depend upon the choice of vector bundle chart and we
note that the form of (\ref{deriv}) automatically ensures that the
covariant derivative satisfies properties (1)--(3) above.

In the special case where the principal bundle $P(M,G)$ is the frame bundle
$LM$ then the structure group $G$ is $GL(n,\Real)$. The tangent bundle is
obtained as an associated bundle simply by taking the canonical
representation of $GL(n, \Real)$ on $\Real^n$ so that $\rho=\hbox{id}$, and
a local trivialisation of $TM$ is given by a system of local coordinates
$x^i$ on $M$. A linear connection on the frame bundle then gives a
covariant derivative defined by 
\begin{equation}
\nabla_XV=\left(X^i\partial_iV^j+\Gamma^j_{ki}V^iX^k\right)\partial_j 
\end{equation}
where the connection coefficients are given by \begin{equation}
\nabla_{\partial_i}\partial_j=\Gamma^k_{ij}\partial_k \end{equation}

We now extend these ideas to the generalised context. A generalised field
$V$ is an element of the $\gs$-module of generalised sections
$\Gamma_\gs(M,E)$. We first extend the definition of generalised linear
connection on a manifold given in \cite{gprg} to the case of a
generalised connection on a vector bundle $E$.  
\begin{definition} \label{vb connection}
(i) A {\it generalised vector bundle connection} $\hat\nabla$ on a vector
bundle $E$ is a map \beas \gs^1_0(M) \times \Gamma_\gs(M,E) &\to&
\Gamma_\gs(M,E) \\ (X,V) &\mapsto& \hat\nabla_XV \\ \eeas satisfying the
following conditions
\begin{itemize}
\item[(1)] $\hat \nabla_XV$ is $\gR$-linear in $V$,

\item[(2)] $\hat \nabla_XV$ is $\gs(M)$-linear in $X$,

\item[(3)] $\hat \nabla_X(\lambda V)=\lambda\hat \nabla_XV+X(\lambda)\cdot
V, \quad \forall \lambda \in \gs(M)$.
\end{itemize}

(ii) Let $(V_\alpha, \Psi_\alpha)$ be a vector bundle chart for $E$ with
coordinates $x^i$, $i=1,\dots, m$ on $M$ and $E_A$ $A=1,\dots n$ the fields
on $E$ induced by the canonical basis on $\Real^n$. We define the {\it
generalised connection coefficients} for this chart to be the $mn^2$
functions $\hat \Gamma^B_{iA} \in \gs(V_\alpha)$ given by \begin{equation} \hat
\nabla_{\partial_i}e_A=\hat \Gamma^B_{iA}e_B, \quad 1 \leq i \leq m, \quad
1 \leq A,B \leq m\ .  \end{equation} 
\end{definition} 
Since $C^\infty(M)$ is a submodule of
$\gs(M)$ and the sheaf $\gs(M)$ is fine (2) and (3) imply the
localisability of any generalised vector bundle connection with respect to
its arguments. (cf. \cite{gprg}).

We now show how a generalised connection on $P(M,G)$ may be used to define
a generalised vector bundle connection $\hat\nabla$ on $E$. The key point
is that because of Theorem 3.2 given a connection $\omega \in
\Omega^1_\gs(P,\g)$ we may (locally) take representatives
$(\omega_\eps)_\eps$ of $\omega$ which are themselves connections.  Given a
bundle chart $(V_\alpha, \Psi_\alpha)$ for $E$ we may as before define the
one parameter family of functions \begin{equation} \nabla_{\eps,
\partial_i}e_A=\Gamma^B_{\eps, iA}e_B \end{equation} It is a straightforward
computation to show using (\ref{coeffs}) that $\hat
\Gamma^B_{Ai}=\left(\Gamma^B_{\eps, iA}\right)_\eps$ defines an element of
$\gs(V_\alpha)$.

If we now define 
\begin{equation} 
\nabla_{\eps, X_\eps}
V_\eps=\left(X_\eps^i\partial_iV^A_{\eps} +\Gamma^A_{\eps,iB}V^B_\eps
X^i_\eps \right)e_A \label{genderiv} \label{veccov}
\end{equation} 
then this defines a generalised
vector bundle connection by setting $\hat\nabla_XV= [\left(\nabla_{\eps,
X_\eps}V_\eps\right)_\eps]$. Furthermore the form of equation
(\ref{genderiv}) guarantees that the generalised covariant derivative
satisfies conditions (1)--(3).

As remarked earlier an important example of a connection on an associated
bundle is provided by a linear connection on $TM$ regarded as an associated
bundle of $LM$. In this case local coordinates $x^i$ on $M$ provide a local
section $s:M \to LM$ by associating with $x \in U \subset M$ the frame 
$\{\partial_i\}_{i=1}^n$. If $\omega_\eps$ is a generalised connection
1-form on $LM$ with values in $gl(n,\R)$  then we define the generalised 
Christoffel symbols $\Gamma^i_{jk,\eps}$ by
\beq
s^*\omega_\eps(V)=\Gamma^i_{jk,\eps}V^jE^k_i \label{genchristoffel}
\eeq
where $X=X^i\partial_i$ and $E^i_j$ is the basis for $gl(n,\R)$ given by
the matrix with a one in the $i$-the column and $j$-th row. Then in
accordance with (\ref{veccov}) this defines a generalised covariant
derivative $\hat \nabla$ by
\begin{equation}
\hat \nabla_{\eps,X}V=\left(X^i\partial_iV^j+\Gamma^j_{\eps,ki}V^iX^k\right)\partial_j 
\end{equation}
which satisfies conditions (1)--(3) of Definition \ref{vb connection}. 
Since $E=TM$ in this case 
these are precisely conditions D1--D3 of Definition 5.1 of 
\cite{gprg}. Hence a generalised connection form on $LM$ defines a
generalised linear connection on $TM$ according to the definition of
\cite{gprg}.

\section{Characteristic Classes for Generalised Connections} \label{characteristic}

The curvature form $\Omega$ of a connection on a principal bundle $P$
together with an invariant polynomial $f$ may be used to construct a closed
2-form $\bar f(\Omega)$ on the base $M$.  It turns out that the element of
the de Rham cohomology $H^*(M)$ associated with $\bar f(\Omega)$ does not
depend upon the choice of connection form $\omega$ on $P$ and is a
topological invariant known as a characteristic class. For the case of a
connection on a vector bundle with fibre $V$ one can view the bundle as an
associated bundle of a principal $GL(V)$ bundle 
and calculate the characteristic classes of this. 
An important example of this are the Chern classes whose
construction we describe below.

We start by showing how to construct closed forms on $M$. Let $f$ be a
$k$-linear symmetric map from $\g$ to $\R$ which is $\ad$-invariant, so that
\begin{equation} 
f(\Ad g V_1,\ldots,\Ad g V_k)=f(V_1,\ldots, V_k), \qquad \forall g \in
G, \quad \forall V_i \in \g \, .  
\end{equation} 
We now define the $2k$-form $f(\Omega)$ on $P$ by 
\begin{equation}
f(\Omega)(v_1,\ldots,v_{2k})=\frac{1}{(2k)!}\sum_\sigma \sign\sigma
f(\Omega(v_{\sigma(1)},v_{\sigma(2)}),\ldots,
\Omega(v_{\sigma(2k-1)},v_{\sigma(2k)}))\ , \label{f} 
\end{equation} 
where $v_i \in T_pP$ and the sum is over all permutations $\sigma$.

We may use $f(\Omega)$ to define a $2k$-form $\bar f(\Omega)$ on $M$
as follows. Let $x \in M$ and $u_i \in T_xM$. Now let $p \in P$ be such
that $\pi(p)=x$ and let $v_i \in T_pP$ be such that $D\pi_p(v_i)=u_i$, then
we define 
\begin{equation} 
\bar f(\Omega)_x(u_1,\ldots, u_{2k})=f(\Omega)_p(v_1,\ldots,v_{2k}) \ .
\label{fbar}
\end{equation} 
We note that the RHS does not depend upon the choice of $v_i$ which
projects onto $u_i$ since any two such vectors differ by a vertical vector
and $\Omega$ vanishes on vertical vectors, nor does the RHS depend upon the
choice of $p$ since any two such points $p_1$ and $p_2$ are related by
$p_2=R_gp_1$ for some $\g \in G$ and $R^*_g\Omega=\Ad(g^{-1})\Omega$. Thus
$f(\Omega)$ projects onto a unique $2k$-form $\bar f(\Omega)$ on $M$.

We next observe that $\bar f(\Omega)$ is closed. This follows from the fact
that $d(f(\Omega))=D(f(\Omega))$ since $\Omega$ vanishes on horizontal
vectors and 
$D(f(\Omega))=0$ by the Bianchi identity. 

Finally using a homotopy argument one can show that given two connection
1-forms on $P$, $\omega_0$ and $\omega_1$ then 
\begin{equation} 
\bar f(\Omega_0)-\bar f(\Omega_1)=dQ \label{Q}
\end{equation} 
where $Q$ depends upon both $\omega_0$ and
$\omega_1$. Hence $\bar f(\Omega_0)$ and $\bar f(\Omega_1)$ represent the
same cohomology class, a so-called characteristic class 
(see \cite{CDD} for details).

The results of section \ref{gencurve} show that the above
constructions may be extended to the case of a generalised connection.

\begin{theorem} \label{charclass} 
Let $\omega \in \Omega^1_{\G}(P,\g)$ be a generalised
connection 1-form on $P$ with generalised curvature $\Omega$ and $f$ an
$\ad$-invariant $k$-multilinear map as above, then $f(\Omega)$ projects onto a
unique closed $2k$-form $\bar f(\Omega) \in \Omega^{2k}_\gs(M)$ which satisfies
\begin{equation} 
\pi^*\bar f(\Omega)=f(\Omega) \ .  
\end{equation} 
\end{theorem}

\pr 
We define $f(\Omega)$ using equation (\ref{f}). Then by Theorem
\ref{horizontal} and Theorem \ref{adinv} $f(\Omega)$ is an invariant
horizontal generalised form so that it projects onto a unique generalised
$2k$-form on $M$ given by equation (\ref{fbar}).
\ep

The theory of generalised exterior calculus on a manifold is outlined in
\cite{ndg} where it is shown that all the classical operations may be
extended to the generalised case. Using these ideas we may define the
concept of the generalised de Rham cohomology of $M$. The following result
clarifies the relationship between generalized and smooth de Rham cohomology.

\begin{proposition}
  For each $p \geqslant 0$, the following isomorphism of real vector spaces holds:
$$
H^p_{\gs}(M) \cong \gR \otimes_{\R} H^p(M)
$$
\end{proposition}
\pr
For clarity, in this proof we denote by $d$ the usual exterior derivative
on smooth forms and by $d'$ the corresponding map on Colombeau forms.
$H^*_\gs$ is calculated through the following fine resolution of the sheaf of
locally constant Colombeau functions:
$$
0 \longrightarrow \ker(d') \stackrel{d'}{\longrightarrow} \Omega_\gs^0(M)
\stackrel{d'}{\longrightarrow} \Omega_\gs^1(M) \stackrel{d'}{\longrightarrow}\dots
$$  
Consider now the vector spaces $\gR\otimes_\R \Omega^p(M)$. Then for $p=0$,
the kernel of $id\otimes d: \gR \otimes_\R C^\infty(M) \to \gR\otimes_\R \Omega^1(M)$
is the set of locally constant Colombeau functions on $M$, i.e., $\ker d'$.
Moreover, for $p\geqslant 1$ the map
$\sum r_i\otimes \omega_i \mapsto \sum r_i \otimes [\omega_i]$
induces a bijection 
\begin{equation}\label{cohomology}
\ker(id\otimes d)/\mbox{Im}(id\otimes d) \to \gR\otimes_\R H^p(M)
\end{equation}
Thus we arrive at another fine resolution of $\ker d'$ of the form
$$
0 \longrightarrow \ker(d') \stackrel{id\otimes d}{\longrightarrow} \gR \otimes_\R C^\infty(M,\R)
\stackrel{id\otimes d}{\longrightarrow}  \gR \otimes_\R \Omega^1(M)
\stackrel{id\otimes d}{\longrightarrow}\dots
$$
which, by the abstract de Rham theorem (\cite{warner}, ch.\ 5), 
calculates the same cohomology as above.
The claim therefore follows from (\ref{cohomology}).
\ep

\begin{definition} 
Let $Z^p_\gs(M)$ be the 
$\gR$-module of closed generalised $p$-forms 
\begin{equation}
Z^p_\gs(M)=\{A \in \Lambda_{\G}^p(M) | dA=0\} \ , 
\end{equation} 
and let $B^p_\gs(M)$ be the $\gR$-module of exact generalised
$p$-forms
\begin{equation} B^p_\gs(M)=\{A \in \Lambda_{\G}^p(M) | \exists B \in
  \Lambda^{p-1}_\gs(M) \mbox{ s.t. }  A=dB\} \ , 
\end{equation} 
then we define the $p$-th generalised de Rham cohomology module $H^p_\gs(M)$ 
to be given by the quotient 
\begin{equation}
H^p_\gs(M)=Z^p_\gs(M)/B^p_\gs(M) 
\end{equation} 
\end{definition}

We now have the following Theorem

\begin{theorem} \label{charclass2} 
Let $\omega \in \Omega^1_{\G}(P,\g)$ be a
generalised connection 1-form on $P$ with generalised curvature $\Omega$ and
$f$ an $\ad$-invariant $k$-multilinear map as above, then $f(\Omega)$ projects
onto a unique generalised closed $2k$-form $\bar f(\Omega)$ which defines
an element of $H_\gs^{2k}(M)$ which does not depend upon $\omega$.  
\end{theorem}

\pr 
By Theorem \ref{charclass} above we know that $\bar f(\Omega)$
defines an element of $Z_\gs^{2k}(M)$. By Theorem \ref{localrep} we may
take a local representative $(\omega_\eps)_\eps$ of 
$\omega$ such that all $\omega_\eps$ are connections. Using this
and the classical result we see that given two different connections
$\omega_0$ and $\omega_1$ then $\bar f(\Omega_0)-\bar f(\Omega_1)$ defines
an element of $B_\gs^{2k}(M)$. Hence $\bar f(\Omega)$ defines an element of 
$H_\gs^{2k}(M)$ which does not depend upon the choice of generalised 
connection.
\ep

Since a (classical) characteristic class defines an element of the de Rham
cohomology $H^*(M)$ which is independent of the connection one may integrate 
it over a cycle to obtain a number which does not depend upon the
connection. In particular one can integrate a characteristic class in
$H^m(M)$ over the whole manifold to obtain a topological invariant. A well
known example of this is the Euler number of an even dimensional manifold
which is the integral of the Euler class $\gamma$ of the frame bundle over
$M$. For a two dimensional manifold the Euler class is proportional to the
curvature two form and one obtains
\begin{equation} 
\chi(M)=\frac{1}{2\pi}\int_M \Omega
\end{equation}    
However if one is looking at a manifold with boundary this formula needs to
be corrected by adding the integral of the geodesic curvature 1-form 
$\kappa_g$ over $\partial M$,
\begin{equation}
\chi(M)=\frac{1}{2\pi}\int_M \Omega +\frac{1}{2\pi}\int_{\partial M}\kappa_g
\label{gauss-bonnet}
\end{equation}
This extra boundary term is an example of a Chern-Simons term which has to
be included when integrating over a manifold with boundary since by
(\ref{Q}) two different connections induce characteristic classes which differ 
by an exact form $dQ$. By including a boundary integral corresponding to
$Q$ one obtains an expression which does not depend upon the connection.

By Theorem \ref{charclass2} it is clear that these ideas may be extended to
the generalised case to give an expression for the Euler number in terms of
the integral of the generalised curvature, however in this case one should
note that the integral of a generalised form is  $\gR$-valued rather than
$\R$-valued. In fact in the generalised case
it is more useful to turn formula (\ref{gauss-bonnet}) 
round to given an expression for
the generalised curvature over a singular region in terms of the Euler
number and a regular integral over the boundary
\begin{equation}
\int_M \Omega = 2\pi\chi(M) -\int_{\partial M}\kappa_g \ .
\label{gauss-Bonnet2}
\end{equation}

This method was used in \cite{clarke} to show that the generalised
curvature of a cone is associated to a delta function.

In the next section a similar idea will be used to relate monopole charges
of generalised connections with  $U(1)$ Chern classes and generalised 
instanton numbers with $SU(2)$ Chern classes. 
We therefore briefly review the general construction of these classes. 
Let $E$ be a complex vector bundle on a real manifold $M$ with typical
fibre $\C^n$. Then such a bundle may be regarded as an associated bundle of
the principal fibre bundle $P(M,GL(n,\C))$ and any connection on $E$ may be
obtained as the induced connection of a suitable connection on $P$ (see
\cite{CDD} for details). 
One may construct invariant polynomials $f_k$ $k=0,\ldots,n$
on $GL(n,\C)$ by writing the characteristic polynomial of the matrix $A$ as
\begin{equation}
\det\left(\lambda I-A/2\pi i\right)=\sum_{k=0}^n f_n(A)\lambda^{n-k}
\end{equation}
Now let $\omega$ be some generalised connection on $P$ then the 
generalised $k$-th Chern class of the vector bundle $E$ is the generalised 
cohomology class of the unique closed
$2k$-form $c_k(\Omega)$ on $M$ which is defined by requiring that
\begin{equation}
\pi^*c_k(\Omega)=f_k(\Omega) \ .
\end{equation}
For a two dimensional manifold the only Chern number one can obtain is given by
\begin{equation}
C_1(E)=\int_M c_1(\Omega) \ .
\end{equation}
However for a four dimensional manifold it is possible to construct two 
Chern numbers
\beas
C_2(E)&=&\int_M c_2(\Omega) \\
C_1^2(E)&=&\int_M c_1(\Omega)\wedge c_1(\Omega) \\
\eeas

\section{Weakly Singular Solutions of Yang-Mills equations} \label{yangmills}

Historically the first example of a weakly singular solution of
the Yang-Mills equations is Dirac's description of the magnetic
monopole \cite{dirac1}. In this paper Dirac describes an
electromagnetic  field with an electromagnetic 4-potential given in
spherical polar coordinates by
\beq
A=\alpha/2(\cos\theta-1)d\phi \label{potential} 
\eeq 
which he
notes is singular at the origin $r=0$ and also along the negative
$z$-axis where $\theta=\pi$. However he claims that if $\alpha=n \in
\N$ the singularity on the axis is physically unimportant. This is
because the integral around a small loop $\gamma$ round the negative
$z$-axis is then given by 
\beq 
\int_\gamma A=2n\pi \label{quant}
\eeq 
so that the singularity has the effect of changing the phase
of the field by an integer multiple of $2\pi$ which does not
effect the physical field. If one now computes the electromagnetic
field $F=dA$ produced by A in this case one finds that 
\beas
F&=&-n/2\sin\theta d\theta\wedge d\phi \\
&=&\frac{n}{2r^3}(xdz\wedge dy+ydx\wedge dz +zdy\wedge dx) \\
\label{emfield} 
\eeas 
which corresponds to an electric and
magnetic field given by 
\beq 
{\bf E}=0, \qquad \hbox{and} \quad {\bf H}=\frac{n{\bf r}}{r^3} 
\eeq 

So that this solution has no
electric field but a magnetic field which looks like an integer
magnetic charge centred at the origin. Indeed if we integrate
over the 2-sphere $S$ given by $t=const.$ and $x^2+y^2+z^2=R^2$
one finds 
\beq 
\int_S {\bf H}.d{\bf S}=4\pi n \label{charge} 
\eeq
On the other hand for a smooth vector field
\beq 
\int_S {\bf H}.d{\bf S} = \int_V \Div {\bf H} dv 
\eeq 
where $V$ is the interior of $S$. Now $\Div(\frac{{\bf
r}}{r^3})=0$ for ${\bf r} \neq 0$, so that $\Div {\bf H}$ behaves
like the Dirac delta function $4\pi n\delta_0$ and the solution is
said to describe a magnetic monopole. As Dirac pointed out the
condition that the integral (\ref{quant}) is an integer multiple
of $2\pi$ leads to the quantisation of magnetic charge according to
(\ref{charge}). However one is not justified in applying
integral theorems in this case due to the singularities of $A$. 

In modern treatments of the magnetic monopole one regards $i{\bf
A}$ as giving a connection on a $U(1)$ bundle on $\R^4 \setminus
\{r=0\} \approxeq \R^2 \times S^2$. One then avoids the singularity
on the negative $z$-axis by giving different descriptions of the
potential on $U_+=\R^4\setminus \{r=0 \hbox{ or }  z<0\}$ and on
$U_-=\R^4\setminus \{r=0 \hbox{ or }  z>0\}$ according to 
\beas
A=A_+&=& n/2(\cos\theta-1)d\phi, \quad\hbox{on}\quad U_+ \\
A=A_-&=& n/2(\cos\theta+1)d\phi, \quad\hbox{on}\quad U_- \\
\eeas 
Note that $iA_+$ and $iA_-$ are related by the gauge transformation 
$iA_+=iA_- -ind\phi$ on the
thin `overlap' $U_+ \cap U_-$ and so define a $U(1)$ connection on
$\R^2\times S^2$ and
that both $A_+$ and $A_-$ agree with the field (\ref{emfield}) on
$U_+ \cup U_-$. In terms of this description the magnetic charge
associated with the monopole is then given in terms of the Chern
class (see \cite{CDD} for details).

However from the spacetime point of view there is nothing wrong at the
origin and it would be preferable to
work with a space with topology $\R^4$ rather than exclude the
origin and work with a space of topology $\R^2\times S^2$. More
importantly the above description relies upon the quantisation
condition $\alpha \in \N$ to ensure that $A_+$ and $A_-$ define a
$U(1)$ connection and this method can not be applied to the sort
of weakly singular Yang-Mills connections we consider later 
which have fractional
charges. We therefore give an alternative formulation
of the Dirac monopole in terms of a
singular Yang-Mills connection over $\R^4$. 

Before doing so we consider the simpler case of a $U(1)$ connection 
which is given by
\beq
A=i\alpha d\phi \label{simple}
\eeq
where $\alpha \in \R$.
It is readily verified that away from the $z$-axis this connection is
flat. However if one considers the holonomy generated by a small loop
$\gamma$ encircling the axis one finds that it is given by
$\exp(-2\pi\alpha i)$ so that unless $\alpha \in \Z$ the loop generates a
non-trivial element of holonomy. Furthermore if we consider the gauge
equivalent connection given by
\beq
\tilde A=A+g^{-1}dg
\eeq
where $g(\phi)=\exp(-in\phi)$ (with $n \in \Z$), 
then $\tilde A$ generates the element of holonomy 
$\exp(-2\pi(\alpha-n))$ and without loss of generality we may use such 
a change of gauge to ensure that $\alpha$ lies in the range $0 \leq \alpha
<1$. If $\alpha=0$ then in this gauge the connection vanishes everywhere 
and is globally flat. 
If $\alpha \neq 0$ the holonomy is non-trivial and the connection cannot be 
extended to a regular connection on the whole of Minkowski space (including the
$z$-axis). This is because in the simple $U(1)$ case the holonomy may be
given in terms of the exponential of the integral of the curvature and so
for a regular connection this tends to the identity as the area of the loop 
shrinks to zero. 

Although it is not possible to give a description of (\ref{simple}) as a
regular connection on the whole of Minkowski space it is possible to define
a generalised connection on the whole of $\R^4$ which represents
$A$. To do this we look at $A$ in Cartesian coordinates. In these
coordinates $A$ is given by
\beq
A=\frac{xdy-ydx}{x^2+y^2}i \alpha
\eeq
This is singular on the $z$-axis so we replace $A \in
\Omega^1(\R^4\setminus \R^2)$ by a regular family of potentials
$A_\eps$ which represent an element in  $\Omega^1_\gs(\R^4,i\R)$
\beq
A_\eps=\frac{xdy-ydx}{x^2+y^2+\eps^2}i \alpha
\eeq
The corresponding field $F_\eps=dA_{\eps}$ is then given by
\beq
F_\eps=\frac{2i\eps^2\alpha dx\wedge dy}{(x^2+y^2+\eps^2)^2}
\eeq
If we now write this in cylindrical polar coordinates $(t,\rho, \phi, z)$
we find
\beq
F_\eps=\frac{2i\eps^2\alpha\rho d\rho\wedge d\phi}{(\rho^2+\eps^2)^2} \label{F}
\eeq
Then by undertaking a calculation very similar to that in \cite{clarke}
one can show that the generalised two form in $\Omega^2_\gs(\R^4,i\R)$
represented  by the family $(F_\eps)_\eps$ given by (\ref{F}) is
associated to  the distributional 2-form 
$2\pi i \alpha \delta^{(2)}(x,y)dx\wedge dy$, which corresponds to a magnetic
field given in Cartesian coordinates by
\beq
{\bf H}=4\pi \alpha (0, 0, \delta^{(2)}(x,y))
\eeq

We now turn to the generalised description of the Dirac monopole. As in the
previous example we start by looking at the Dirac potential $A$ given by
(\ref{potential})  in Cartesian coordinates. In these
coordinates $A$ is given by
\begin{equation}
A=\frac{\alpha}{2}
\left(\frac{z}{(x^2+y^2+z^2)^{1/2}}-1\right)
\left(\frac{xdy-ydx}{x^2+y^2}\right)
\end{equation}
We first note that in these coordinates the coefficients of
$dx$ and $dy$ both diverge as we approach the negative $z$-axis.
To remedy this we again replace $A$ by a regular family of potentials
$(A_\eps)_\eps $ which represent an element of $\Omega^1_\gs(\R^4)$
\beq
A_\eps=\frac{\alpha}{2}\left(\frac{z}{(x^2+y^2+z^2+\eps^2)^{1/2}}-1\right)
\left(\frac{xdy-ydx}{x^2+y^2+\eps^2}\right)
\eeq 
The associated electromagnetic field is given by
$F_\eps=dA_\eps$ and a direct calculation shows that 
\beq 
F_\eps=F^1_\eps+F^2_\eps+F^3_\eps
\eeq
where
\beas
F^1_\eps&=&\frac{\alpha (xdz\wedge dy+ydx\wedge dz +zdy\wedge dx)}
{2(r^2+\eps^2)^{3/2}} \\
F^2_\eps&=&\frac{\alpha\eps^2 zdx\wedge dy}
{2(r^2+\eps^2)^{3/2}(x^2+y^2+\eps^2)} \\
F^3_\eps&=& \left(\frac{z}{(r^2+\eps^2)^{1/2}} - 1 \right)
\frac{\alpha\eps^2 dx\wedge dy}{(x^2+y^2+\eps^2)^2}\\ 
\eeas
The first term $F^1_\eps$ gives the monopole expression (\ref{emfield}) as
$\eps \to 0$ while the flux of the second term vanishes as $\eps \to 0$. 
However the
third term $F^3_\eps$ diverges on the negative $z$-axis and its
contribution must be taken into account in computing the flux integral.
Now away from the negative $z$-axis $F^3_\eps \to 0$ as $\eps \to 0$, so
that
\beas
\lim_{\eps \to 0}\int_S F^3_\eps&=&\lim_{\eps \to 0}\int_{S^+} F^3_\eps +
\lim_{\eps \to 0}\int_{S^-} F^3_\eps \\
&=&0-\lim_{\eps \to 0}\int_{S^-}\frac{2\alpha\eps^2 dx\wedge dy}
{(x^2+y^2+\eps^2)^2} \\
&& +\lim_{\eps \to 0}\int_{S^-}\left(1+\frac{z}
{(r^2+\eps^2)^{1/2}}\right)
\frac{\alpha\eps^2 dx\wedge dy}{(x^2+y^2+\eps^2)^2} \\
&=&0-\lim_{\eps \to
  0}\oint_C\alpha\left(\frac{xdy-ydx}{x^2+y^2+\eps^2}\right) +0 \\
&=&-\int_0^{2\pi}\alpha d\phi \\
&=& -2\alpha\pi \\
\eeas

We therefore see that in the regularised description as well as the
monopole term $F^1$ there is also another term $F^3$ which behaves like the
flux due to a potential $\alpha d\phi$ which corresponds to a delta function 
like magnetic field along the negative $z$-axis. It is this term due to the
`wire singularity' that results in the Chern number vanishing 
(as it must do given the local
trivialisation over $\R^4$) despite the presence of the monopole term.
In terms of a $U(1)$ bundle description, if $\alpha=n\in \Z$,
the connection $ind\phi$ generates a trivial holonomy $e^{2m\pi i}=1$ 
when integrated round any closed curve $\gamma$ and it is this feature 
that enables one to give the usual modern description of
a monopole in terms of a $U(1)$ bundle over $\R^2 \times S^2$. 

A more interesting example of a weakly singular connection is the
``fractionally charged instanton'' discovered by Forgacs et al
\cite{Forgasc1} (see also \cite{Forgasc2}, \cite{Chak}, \cite{Braam}). This
is a (Euclidean) self-dual $SU(N)$ Yang-Mills connection on the 4-sphere with a
singularity along a 2-sphere. Rather than consider the exact instanton
solutions we follow R\aa de \cite{Rade} and look at the general class of weakly
singular $SU(N)$ connections on the unit ball in Euclidean
space. Furthermore  for simplicity of exposition we will look in detail at
the case $N=2$ since no new features arise for higher values of $N$. 

Let $B$ denote the 4-ball $x^2+y^2+z^2+w^2 \leq 1$ in Euclidean
4-space with $(x,y,z,w)$ Cartesian coordinates. Let $D$ denote the disk
$z=0$, $w=0$, $x^2+y^2 \leq 1$. 
We may now define 4-dimensional cylindrical polar 
coordinates $(x,y,r,\phi)$ by $z=r\cos\phi$ and $w=r\sin\phi$. Let $P$ be
an $SU(2)$ bundle over $B$, then by a
singular $SU(2)$ connection on $P$ we mean a connection that in a local
gauge may be written 
\beq
\omega=\tilde \omega+a \label{singconn}
\eeq
where
\beq
\tilde \omega=\left(\begin{array}{cc}
i\alpha & 0 \\
0 & -i\alpha \\
\end{array}\right)d\phi, \qquad 2\alpha \not\in \Z,
\label{coneconn}
\eeq      
and $a$ is a regular $su(2)$-valued 1-form.  
The first term defines a regular connection on $B \setminus D$, but the
holonomy around a loop  $\gamma$ round the disk $D$ is given by
\beq
\left(\begin{array}{cc}
\exp(-2\pi\alpha i) & 0 \\
0 & \exp(2\pi\alpha i) \\
\end{array}\right)
\eeq
and it is the non-trivial nature of this holonomy that makes the connection
singular on $D$. (Note the reason that we demand that $2\alpha$ rather than
$\alpha$ is not an integer is that in the $SU(2)$ case it is possible to
regard $A$ as splitting as the direct sum of a trivial $SU(2)$ 
connection on $B$ and a flat $U(1)$ connection $2i\alpha d\theta$ on 
$B \setminus D$ with holonomy $\exp(-4\pi\alpha i)$ round $D$, see
\cite{Rade} for details).

More generally one can consider singular connections on 4-manifolds which
have non-trivial limit holonomy along an embedded surface in the
4-manifold. A theorem by L Sibner and R Sibner \cite{sibnersibner} shows
that a finite energy Yang-Mills connection on the complement of an embedded
surface in a 4-manifold has a well defined limit holonomy around the
surface and that the connection can be extended across the surface if and
only if the limit holonomy is trivial. Furthermore they also show that 
any such connection over
$B \setminus D$ that is in $L^{2,1}$ locally and has curvature in $L^p$
globally is gauge equivalent to (\ref{singconn}) (see also \cite{Rade2}) so
there is no loss of generality in considering singular connections given by
(\ref{singconn}).   

We now show how to replace (\ref{coneconn}) by a generalised connection on
$B$ and hence also give a description of (\ref{singconn}) as a generalised
connection on $B$. This enables us to give an explicit formula for the 
curvature of the regularised version of (\ref{singconn}). 
As usual we start by writing (\ref{coneconn}) in Cartesian coordinates  
\beq
\tilde \omega=\left(\begin{array}{cc}
i\alpha & 0 \\
0 & -i\alpha \\
\end{array}\right)
\frac{xdy-ydx}{x^2+y^2}
\label{coneconnb}
\eeq      
and introduce the corresponding regular family of connections according to  
\beq
\tilde\omega_\eps=\left(\begin{array}{cc}
i\alpha & 0 \\
0 & -i\alpha \\
\end{array}\right)
A_\eps
\label{regconeconnb}
\eeq
where $A_\eps$ is the 1-form given by
\beq
A_\eps=\frac{xdy-ydx}{x^2+y^2+\eps^2} \label{newA}
\eeq
This may be used to define a generalised connection in
$\Omega^1_\gs(P,su(2))$ which is represented by the family
$(\omega_\eps)_\eps$ given in a local gauge by 
$\omega_\eps=\tilde\omega_\eps+a$.
The corresponding curvature is an element of $\Omega^2_\gs(P,su(2))$ which
is represented by the family given by 
\beq
F_\eps=d\tilde\omega_\eps+da+[\tilde\omega_\eps+a,  \tilde\omega_\eps+a]
\eeq
In order to calculate this we decompose the $su(2)$ valued 1-form $a$ into
its diagonal and off diagonal components 
\beq
\left(\begin{array}{cc}
a_D & a_T \\
-\bar a_T & -a_D \\
\end{array}\right)
\eeq
and use the fact that if $a$ and $b$ are two $su(2)$-valued 1-forms then
\beas
{[a,b]}_D&=&-2i\hbox{Im}(a_T\wedge\bar b_T) \\
{[a,b]}_T&=&2(a_D\wedge b_T-a_T\wedge b_D) \\
\eeas
We then find that
\beq
F_\eps=F^1_\eps+F^2
\eeq
where the first term is the singular part and is 
given in terms of its diagonal and off diagonal pieces by
\beq
F^1_{\eps,D}=dA_\eps, \qquad F^1_{\eps, T}=-2i\alpha a_T \wedge A_\eps 
\label{susing}
\eeq 
with $A_\eps$ given by (\ref{newA}) and $dA_\eps$ given by
\beq
dA_\eps=\frac{2\eps^2 dx\wedge dy}{(x^2+y^2+\eps^2)^2}
\eeq
while the second term is the regular part and is just the curvature of the
regular part of the connection
\beq
F^2=da+[a,a] \label{smooth}
\eeq
The curvature of the generalised connection therefore splits into a smooth
part $F^2 \in \Omega^2(P,su(2))$ given by (\ref{smooth}) and a singular
part $F^1 \in \Omega^2_\gs(P,su(2))$ represented by the family given by 
(\ref{susing}), the new feature in the $su(2)$ case being the way that the
off diagonal piece is the wedge product of the smooth and singular parts of
the connection.  

{\it Acknowledgment:} We would like to thank Andreas Cap for helpful discussions.


\end{document}